\documentclass[a4paper,12pt]{article}
\usepackage{amssymb}
\usepackage{amsthm}
\newtheorem{Lemma}{Lemma}
\newtheorem{Proposition}{Proposition}
\newtheorem{Theorem}{Theorem}

\newtheorem{Definition}{Definition}
\newtheorem{Claim}{Claim}

\newcommand{\PP}{\mathbb{P}}
\newcommand{\C}{\mathbb{C}}

\begin{document}
\title{On a Lemma of Bompiani}
\author{E.~Ballico and C.~Fontanari\thanks{This research is part of the 
T.A.S.C.A. project of I.N.d.A.M., supported by P.A.T. (Trento) and 
M.I.U.R. (Italy).} }
\date{}
\maketitle 

\begin{small}
\begin{center}
\textbf{Abstract}
\end{center}
We reprove in modern terms and extend to arbitrary dimension a classical 
result of Enrico Bompiani on algebraic surfaces $X \subset \PP^r$ having 
very degenerate osculating spaces.
\end{small}
\vspace{0.5cm}

Let $X \subset \PP^r$ be an integral nondegenerate projective 
variety of dimension $n$ defined over the field $\C$. 
In order to introduce the osculating space of order $m$ to a 
smooth point $p \in X$, fix a lifting 
\begin{eqnarray*}
U \subseteq \C^n &\longrightarrow& \C^{r+1} \setminus \{ 0 \} \\
t &\longmapsto& p(t)
\end{eqnarray*}
of a local parametrization of $X$ centered in $p$ and define
$T(m,p,X)$ to be the linear span of the points $[p_I(0)] \in \PP^r$, 
where $I$ is a multi-index such that $\vert I \vert \le m$. 
The starting point of our research was the following result, which 
we read in a recent paper by Luca Chiantini and Ciro Ciliberto (see 
\cite{ChiCil:02}, Proposition~2.3): 

\begin{Proposition}
Let $X \subset \PP^r$ be a smooth variety and let $p \in X$ be a general 
point. Assume that $\dim T(m,p,X) = \dim T(m+1,p,X) = h$. 
Then $X \subseteq \PP^h$. 
\end{Proposition}

Our natural question was: can one go a little bit further? 
Namely, if the dimension of the osculating space at a general point 
does not increase too much while passing from order $m$ to order $m+1$, 
what can one say about the projective geometry of $X$? 
Following a suggestion of Ciro Ciliberto, to whom we are grateful, 
we looked for an answer among the works of Bompiani (see 
\cite{Bompiani:78}). Enrico Bompiani (1889--1975) was a student of 
Guido Castelnuovo and in more than three hundreds papers intensively 
studied the differential geometry of projective varieties; in 
particular, he deeply investigated the relationship between partial  
differential equations and algebraic geometry. In the beautiful paper 
\cite{Bompiani:19} we found an explicit answer to our question. 
Unluckily (or fortunately, according to the points of view), 
Bompiani's explanation turns out to have a couple of drawbacks: 
first, he treats only the case of surfaces; and next, his arguments 
are a little bit involved and it is not so easy to understand them 
properly. Therefore we decided to work them out again in a hopefully 
more clear and rigorous form; as a result, we obtained the following 
generalization of Bompiani's Lemma (see \cite{Bompiani:19}, pp.~614--615):

\begin{Theorem}~\label{bompiani}
Let $X \subset \PP^r$ be a smooth variety and let $p \in X$ be a general 
point. Assume that $\dim T(m,p,X) = h$ and $\dim T(m+1,p,X) = h + k$ with 
$1 \le k \le n-1$. 
Then either $X \subset \PP^{h+k}$ or $X$ is covered by infinitely 
many subvarieties $Y$ of dimension at least $n-k$ such that 
$Y \subset \PP^{h-m}$.
\end{Theorem}

Theorem~\ref{bompiani} is a direct consequence of the following 

\begin{Claim}\label{main}
Under the assumptions of Theorem~\ref{bompiani}, either 
$X \subset \PP^{h+k}$, or there exists a subvariety $Y$ of $X$ 
such that $p \in Y$, $\dim Y \ge n-k$, and $T(m,q,X) = T(m,p,X)$ 
for a general $q \in Y$. 
\end{Claim}

Indeed, assume for a moment that Claim~\ref{main} holds. 
Since $Y \subset T(m,p,X)$, we have $<(m+1)Y> \subseteq 
T(m,p,X)$. Moreover, from \cite{ChiCil:02}, Remark~2.4, it 
follows that either $<(m+1)Y> = \PP^r$ or 
$\dim <(m+1)Y> \ge \dim <Y> + m$. 
Hence we deduce that either 
$$
X \subset \PP^r = <(m+1)Y> \subseteq \PP^{h},
$$ 
or  
$$
\dim <Y> \le \dim <(m+1)Y> - m \le h - m. 
$$    
Therefore we are reduced to establish the claim. In order to do that, 
we recall a standard definition:

\begin{Definition}
The order $m$ osculating variety of $X$ is 
$$
T(m,X) := \overline{\bigcup_{p \in X} T(m,p,X)}.
$$
\end{Definition}

We also recall a natural description of the tangent space
to $T(m,X)$ at a general point $P \in T(m,p,X)$: 
$$
T_P(T(m,X)) = \left\{ \frac{d \gamma}{d s}(0) \in \PP^r \right\}
$$
where 
\begin{eqnarray*}
\gamma: \Delta &\longrightarrow& T(m,X) \subset \PP^r \\
s &\longmapsto& \sum_{\vert I \vert \le m} \alpha_I(s) [ p_I (t(s))]
\end{eqnarray*}
is a holomorhic map defined on the unit disc $\Delta \subset \C$
(just notice that $P$ is a smooth point of the integral variety 
$T(m,X)$ and see \cite{GriHar:78}, p.~22). 
We are going to deduce Claim~\ref{main} from the following Lemma, 
which we believe to be quite interesting also in its own:  

\begin{Lemma}\label{include}
Let $X$ be a smooth variety and let $P \in T(m,p,X)$ be a general point
of $T(m,X)$. Then
$$
T_P(T(m,X)) \subseteq T(m+1,p,X).
$$ 
\end{Lemma}

\proof Indeed, we have
\begin{eqnarray*}
\frac{d \gamma}{d s}(0) &=& \sum_{\vert I \vert \le m} 
\dot \alpha_I(0) [ p_I (0)] + \sum_{\vert I \vert \le m} 
\alpha_I(0) \sum_{j=1}^n [ \frac{\partial}{\partial t_j} p_I (0)]
\dot t(0) = \\
&=& \sum_{\vert I \vert \le m+1} \beta_I [ p_I (0)]
\in T(m+1,p,X)
\end{eqnarray*}
and the proof is over.

\qed

Under the assumptions of Theorem~\ref{bompiani}, we have 
$$
\dim T(m,X) \le \dim T(m,p,X) + k
$$ 
by Lemma~\ref{include}. 
It follows that the tangent space $T(m,p,X)$ is constant along a 
subvariety of dimension at least $n-k$, hence Claim~\ref{main} 
holds.

\vspace{0.5cm}


\vspace{0.5cm}

\noindent
Edoardo BALLICO \newline
Universit\`a degli Studi di Trento \newline
Dipartimento di Matematica \newline
Via Sommarive 14 \newline
38050 Povo (Trento) \newline
ITALIA \newline
e-mail: ballico@science.unitn.it

\vspace{0.5cm}

\noindent
Claudio FONTANARI \newline
Universit\`a degli Studi di Trento \newline
Dipartimento di Matematica \newline
Via Sommarive 14 \newline
38050 Povo (Trento) \newline
ITALIA \newline
e-mail: fontanar@science.unitn.it

\end{document}